\documentclass[12pt]{article}
\usepackage{amssymb,color}
\usepackage{latexsym}
\usepackage{amsfonts}
\usepackage{amsmath}

\newtheorem{theorem}{Theorem}[section]
\newtheorem{lemma}[theorem]{Lemma}
\newtheorem{proposition}[theorem]{Proposition}

\newtheorem{corollary}[theorem]{Corollary}
\newtheorem{remark}[theorem]{Remark}

\newtheorem{example}[theorem]{Example}
\numberwithin{equation}{section}

\def\P{\mathbb{P} }
\def\E{\mathbb{E} }
\def\R{\mathbb{R} }
\def\N{\mathcal{N}}
\def\mE{{\rm E}}
\def\mP{{\rm P}}
\def\R{\mathbb{R}}
\def\d{{\rm d}}

\begin{document}

	\allowdisplaybreaks
	\title{ Weak convergence of the extremes of  branching L\'evy processes with regularly varying tails}
\author{ \bf  Yan-Xia Ren\footnote{The research of this author is supported by the National Key R\&D Program of China (No. 2020YFA0712900) and NSFC (Grant Nos. 12071011  and 11731009).\hspace{1mm} } \hspace{1mm}\hspace{1mm}
Renming Song\thanks{Research supported in part by a grant from the Simons
Foundation (\#960480, Renming Song).} \hspace{1mm}\hspace{1mm} and \hspace{1mm}\hspace{1mm}
Rui Zhang\footnote{
Corresponding author.
The research of this author is supported by NSFC (Grant No. 11601354), Beijing Municipal Natural Science Foundation(Grant No. 1202004), and Academy for Multidisciplinary Studies, Capital Normal University.}
\hspace{1mm} }
\date{}
\maketitle
	
\begin{abstract}
In this paper, we study the weak convergence of the extremes of supercritical branching L\'evy processes $\{\mathbb{X}_t, t \ge0\}$
whose spatial motions are L\'evy processes with regularly varying tails.
The result is drastically different from
the case of branching Brownian motions. We prove
that, when properly renormalized, $\mathbb{X}_t$ converges weakly.
As a consequence, we obtain a limit theorem for the order statistics
of $\mathbb{X}_t$.
\end{abstract}
\medskip
\noindent {\bf AMS Subject Classifications (2020)}:
Primary 60J80, 60F05;
Secondary  60G57, 60G70

\medskip

\noindent{\bf Keywords and Phrases}: Branching L\'evy process,  extremal process,
regularly varying,  rightmost position.
\begin{doublespace}
	
\section{Introduction}
We consider a supercritical branching L\'evy process. At time $0$, we start with a single particle which
moves according to a L\'evy process
$\{\xi_t,\mP_x\}$ with L\'evy exponent $\psi(\theta)=\log \mE(e^{i\theta \xi_1})$.
The lifetime of each particle is  exponentially distributed with parameter $\beta$, then it  splits into $k$ new particles with probability  $p_k$, $k\ge0$.
Once born, each particle will independently move (according to the same L\'evy process) and split (according to the same offspring distribution).
We use $\P_x$ to denote the law of the branching L\'evy process when the initial particle starts at position $x$.
The expectation with respect to
$\P_x$ and $\mP_x$ will be denoted by $\E_x$ and $\mE_x.$, respectively.
 We write $\P:=\P_0$, $\E:=\E_0$,  $\mP:=\mP_0$ and $\mE:=\mE_0.$

In this paper, we use ``:=" as a way of definition.
For $a, b\in\R$, $a \wedge b := \min\{a, b\}$.
We will label each particle using the classical Ulam-Harris system.  We write $\mathbb{T}$ for the set of all the particles in the tree, $o$ for the root of the tree. For each particle $u$, we introduce some notation.
\begin{itemize}
	\item $b_u$ and $\sigma_u$:  the birth time and death time of $u$ respectively.
	\item  $\{\xi_t^u:t\in[b_u,\sigma_u]\}$:
	the spatial trajectory of $u$.
	\item
	   $\tau_u:=\sigma_u-b_u$ is the life length of $u$ and $\tau_{u,t}:=\sigma_u\wedge t-b_u\wedge  t$ is the life length of $u$ between $[0,t]$.
	\item $\mathcal{F}_t^{\mathbb{T}}:=
		 \sigma\{b_u\wedge t,\sigma_u\wedge t: u\in\mathbb{T}\}.$
	\item  $X_u:=\xi^u_{\sigma_u}-\xi^u_{b_u}$ and $X_{u,t}:=\xi^u_{\sigma_u\wedge t}-\xi^u_{b_u\wedge t}.$
Note that given $\mathcal{F}_t^\mathbb{T}$,
	$X_{u,t}, u\in\mathbb{T},$ are independent, and
	$$X_{u,t}\overset{d}{=}\xi_{\tau_{u,t}}.$$
	\item
	$I_v$: the set of all the ancestors of $v$,
	including $v$ itself.
	\item
	$n_t^v:$ the number  of particles in $I_v\setminus \{o\}$.
	\item
	$\mathcal{L}_t$  is the set of all  particles alive at time $t$ and $Z_t$  is the number of  particles alive at time $t$.
	\end{itemize}
For $t\geq 0$, define $\mathbb{X}_t:=\sum_{u\in\mathcal{L}_t}\delta_{\xi^u_t}$. The measure-valued process $\{\mathbb{X}_t, t\geq 0\}$ is called
a branching L\'evy process.

It is well known that $\{Z_t; t\geq 0\}$ is a continuous time branching process.  In this paper, we consider the supercritical case, that is, $m:=\sum_k kp_k>1$. Then $\P(\mathcal{S})>0,$ where $\mathcal{S}$ is the event of survival.
The extinction probability $\P(\mathcal{S}^c)$ is the smallest root in $(0,1)$ of the equation $\sum_{k}p_k s^k =s$, see, for instance, \cite[Section III. 4]{Athreya-Ney}.
The family $\{e^{-\lambda t}Z_t, t\ge0\}$, where $\lambda=\beta(m-1)$, is a non-negative martingale and hence
$$\lim_{t\to\infty}e^{-\lambda t}Z_t
=:W \quad \mbox{exists a.s.}$$

For any two
functions $f$ and $g$ on $[0,\infty)$, $f\sim g$ as $s\to 0_+$
 means that $\lim_{s\downarrow 0} \frac{f(s)}{g(s)}=1.$
 Similarly,  $f\sim g$ as $s\to \infty$
 means that $\lim_{s\to\infty} \frac{f(s)}{g(s)}=1.$
Throughout this paper we assume the following
two conditions hold. The first condition is on the offspring distribution:
\begin{itemize}
\item[{\bf (H1)}]\quad
$\sum_{k\ge 1} (k\log k)  p_k<\infty.$
\end{itemize}
Condition (H1) ensures that
 $W$ is non-degenerate with $\P(W>0)=\P(\mathcal{S})$.
For more details, see \cite[Section III.7]{Athreya-Ney}.
The second condition is on the spatial motion:
\begin{itemize}
\item[{\bf (H2)}]
There exist a complex constant $c_*$ with $\Re(c_*)>0$, $\alpha\in(0,2)$ and a function $L(x):\R_+\to \R_+$
slowly varying at $\infty$ such that
$$
\psi(\theta)\sim -c_* \theta^{\alpha}L(\theta^{-1}), \qquad \theta\to 0_+.
$$
\end{itemize}
Since $e^{\psi(\theta)}=\mE(e^{i\theta \xi_1})$,
we have
 $\Re(\psi)\le 0$
 and  $\psi(-\theta)=\overline{\psi(\theta)}$. Thus
$$\psi(\theta)\sim -\overline{c_*} |\theta|^{\alpha}L(|\theta|^{-1}), \quad \theta\to 0_-.$$
Under condition  (H2), one can prove that (see Remark \ref{rek:tail} below)
\begin{equation*}
 \mP (|\xi_s|\ge x)\sim  c s x^{-\alpha}L(x),\quad x\to\infty,
\end{equation*}
that is, $|\xi_s|$ has regularly varying tails.

An important example
satisfying (H2) is the strictly stable process.

\begin{example} {\bf (Stable process.)} Let $\xi$
be a strictly $\alpha$-stable process, $\alpha\in (0, 2)$,  on $\R$ with L\'evy measure
$$n(dy)=c_1x^{-(1+\alpha)}{\bf 1}_{(0,\infty)}(x){\rm d} x+c_2|x|^{-(1+\alpha)}{\bf 1}_{(-\infty,0)}(x){\rm d} x,$$
where $c_1,c_2\ge 0$, $c_1+c_2>0,$
and if $\alpha=1,$ $c_1=c_2=c$.
For $\alpha\in(1,2)$, by \cite[Lemma 14.11, (14.19)]{Sato}
and the fact $\Gamma(-\alpha)=-\alpha\Gamma(1-\alpha)$, we obtain that, for $\theta>0$,
\begin{align*}
	\int_0^\infty (e^{i\theta y}-1-i\theta y) n(\d y)=-c_1\alpha\Gamma(1-\alpha) e^{-i\pi\alpha /2}\theta^\alpha,
\end{align*}
and taking conjugate on both sides of \cite[Lemma 14.11 (14.19)]{Sato}, we have that
\begin{equation*}
		\int_{-\infty}^0(e^{i\theta y}-1-i\theta y) n(\d y)=-c_2\alpha\Gamma(1-\alpha)e^{i\pi\alpha /2}\theta^\alpha.
\end{equation*}
Thus the L\'evy exponent of $\xi$ is given
by: for $\theta>0$,
\begin{align}\label{stable-exponent1}
	\psi(\theta)=\int (e^{i\theta y}-1-i\theta y) n(\d y)=
		-\alpha\Gamma(1-\alpha)(c_1e^{-i\pi\alpha/2}+c_2e^{i\pi\alpha/2})
		\theta^\alpha.
\end{align}
Similarly,
by \cite[Lemma 14.11 (14.18),(114.20)]{Sato}, we have for $\theta>0$,
\begin{align}\label{stable-exponent}
\psi(\theta)&=\left\{\begin{array}{ll}
	\displaystyle\int (e^{i\theta y}-1) n(\d y)&\alpha\in (0,1);\\
	\displaystyle \int (e^{i\theta y}-1-i\theta y{\bf 1}_{|y|\le 1}) n(\d y)+i a\theta, &\alpha=1
\end{array}\right.\\
&=\left\{\begin{array}{ll}
\displaystyle-\alpha\Gamma(1-\alpha)(c_1e^{-i\pi\alpha/2}+c_2e^{i\pi\alpha/2})\theta^\alpha,  &\alpha\in(0,1);\\
\displaystyle-c\pi\theta+ia\theta, &\alpha=1,
\end{array}\right.
\end{align}
where $a\in\R$ is a constant.
It is clear that $\psi$ satisfies (H2). For more details about the stable processes, we refer
the readers to \cite[Section 14]{Sato}.
\end{example}

In Section \ref{exam},
we will give  more  examples satisfying condition (H2).
Note that the non-symmetric 1-stable process does not satisfy (H2).
However, in Example \ref{1-stable}, we will show that our main result still holds for the non-symmetric 1-stable process.

The maximal  position $M_t$ of
a branching Brownian motion
has been studied intensively. Assume that  $\beta=1$, $p_0=0$ and $m=2$. In the seminal paper \cite{KPP}, Kolmogorov,  Petrovskii and  Piskounov proved that
$M_t/t\to \sqrt{2}$ in probability as $t\to\infty$.
Bramson  proved in \cite{Bramson78} (see also \cite{Bramson}) that, under some moment conditions,
$\P(M_t-m(t)\le x)\to 1-w(x)$
as $t\to\infty$ for all $x\in \R$, where $m(t)=\sqrt{2}t-\frac{3}{2\sqrt{2}}\log t$ and  $w(x)$ is a traveling wave solution.
For more works on $M_t$, see \cite{Chauvin88,Chauvin,LS87,Robert}.
 For  inhomogeneous branching Brownian motions, many papers discussed the growth rate of the maximal position, see Bocharov and Harris \cite{Bocharov-Harris14,Bocharov-Harris16} and Bocharov \cite{Bocharov} for the case with catalytic branching at the origin, Shiozawa \cite{Shiozawa}, Nishimori et al. \cite{Nishimori-Shiozawa}, Lalley and Sellke \cite{LS88,LS89} for the case with some general branching mechanisms.

Recently, the full statistics of the extremal
configuration of branching Brownian motion were studied.
Arguin et al. \cite{ABK12,ABK} studied the limit property of the extremal process of branching Brownian motion. They proved that the random measure defined by
$$
\mathcal{E}_t:=\sum_{u\in\mathcal{L}_t}\delta_{\xi^u_t-m(t)}
$$
converges weakly, and the limiting process is a (randomly shifted) Poisson cluster process.  Almost at the same time,
A\"{i}d\'ekon et al. \cite{ABBS} proved similar results using a totally different method.

 For  branching random walks,  several authors have studied
 similar problems under an exponential
 moment assumption on the  displacements of the offspring from the parent,
see A\"{i}d\'{e}kon \cite{Aldekon},
Carmona and Hu \cite{Carmona-Hu},
Hu and Shi \cite{HS},
and Madaule \cite{Madaule}.
When  the displacements of the offspring from the parents are i.i.d.  and
have regularly varying tails,
Durrett \cite{Durrett83} studied
the  limit property of its  maximum displacement $M_n$.    More precisely,  Durrett proved that $a_n^{-1}M_n$  converges weakly,  where $a_n=m^{n/\alpha}L_0(m^n)$ and $L_0$ is slowly varying at $\infty$. Recently,  the extremal processes of the branching random walks with regularly varying steps were studied by Bhattacharya et al. \cite{BHR,BHR2}. In \cite{BHR,BHR2},
it was proved
that the point random measures $\sum_{|v|=n}\delta_{a_n^{-1}S_v}$, where $S_v$ is the position of $v$, converges weakly to a Cox cluster process, which are quite different from the case with exponential moments.   See also \cite{BMPR,Gantert} for related works on branching random walks with heavy-tailed displacements.

Shiozawa \cite{Shiozawa2} studied branching symmetric stable processes with  branching rate $\mu$ being a measure on $\R$ in a Kato class and offspring distribution $\{p_n(x), n\geq 0\}$ being
spatially dependent.
Under some conditions on  $\mu$ and $\{p_n(x), n\geq 0\}$, Shiozawa \cite{Shiozawa2} proved that the growth rate of the  maximal displacement is exponential with rate given by the  principal eigenvalue of the mean semigroup of the  branching symmetric stable processes.
In this paper, we study the extremes of branching L\'evy processes
with  constant branching rate $\beta$ and spatial motion having regularly varying tails (see condition (H2)).  Since our branching rate $\beta$ is not compactly supported, one can not get the growth rate of the  maximal displacement from Shiozawa \cite{Shiozawa2}. As a corollary of our extreme limit result, we  get the growth rate of the   maximal displacement, see Corollary \ref{cor:rightmost*} below.

 The key idea of the proof  in this paper is the
 ``one large jump principle"
which was inspired by \cite{BHR, BHR2, Durrett83}.
Along the  discrete times $n\delta$, the branching L\'evy processes $\{\mathbb{X}_{n\delta},n\ge 1\}$ is a
branching random walk and the displacements from parents has the same law as  $\mathbb{X}_\delta$.
It is natural to think that one may get
the results of this paper from the results for branching random walks directly.
However we can not apply the results for branching random walks in \cite{BHR, BHR2,Madaule} to $\{\mathbb{X}_{n\delta},n\ge 1\}$.
First, under condition (H2),  the exponential moment assumption in \cite{Madaule} is not satisfied.
Secondly,  \cite{BHR} assumes that  the displacements
are i.i.d., while the atoms of the random measure $\mathbb{X}_\delta$ are not independent.
Lastly, although the displacements of offspring coming from the same parent are allowed to be dependent in \cite{BHR2},
Assumption 2.5 in  \cite{BHR2}, where the displacements from parents are given by a special form (see \cite[(2.9) and (2.10)]{BHR2}), seems to be very difficult to check for $\mathbb{X}_\delta$.

 Branching L\'evy processes are closely related to the Fisher-KPP equation when
the classical Laplacian $\Delta$
 is replaced by the infinitesimal generator of the corresponding L\'evy process.  For any $g\in C_b^+(\R)$,  define
$u_g(t,x)=\E_{x}\left(e^{-\int g(y)\mathbb{X}_t(\d y)}\right)$.
 By the Markov property and
 branching property,
 we have that
 $$u_g(t,x)=\mE_x\left(e^{-g(\xi_t)}\right)+\mE_x\int_0^t \varphi(u_g(t-s,\xi_s))\mbox{d}s,$$
 where $\varphi(s)=\beta\left(\sum_k s^kp_k-s\right)$. Then $1-u_g$ is a mild solution to
 \begin{equation}\label{F-KPP}
 	\partial_t u-\mathcal{A} u=-\varphi(1-u),
 \end{equation}
 with initial data $u(0,x)=1-e^{-g(x)}$, where $\mathcal{A}$ is the infinitesimal generator of $\xi$. In \cite{CR},
Cabr\'{e} and Roquejoffre proved that,
under the assumption that the density of $\xi$ is comparable to that of a symmetric $\alpha$-stable process,
the front position of $1-u$ is exponential in time.
Using our main result, we  give
another proof of \cite[Theorem 1.5]{CR} and also partially generalize it,
see Remark \ref{speed on solution}.

\subsection{Main results}
Put ${\R}_0=(-\infty,\infty)\setminus \{0\}$, and $\overline{\R}_0=[-\infty,\infty]\setminus\{0\}$.
Let  $C_{b}^0(\R)$ be the set of all  bounded continuous functions vanishing in a neighborhood of $0$.
Let $C_c^+(\overline{\R}_0)$ be the set of
all non-negative continuous functions on $\overline{\R}_0$
such that $g=0$  on $(-\delta,0)\cup(0,\delta)$ for some $\delta>0$.
It is clear that  if $g\in C_c^+(\overline{\R}_0),$ then  $g^*(x):={\bf 1}_{\R_{0}}(x) g(x) \in C_b^0(\R)$.
Denote by $\mathcal{M}(\overline{\R}_0) $ the set of all Radon
measures endowed with the topology of vague convergence (denoted by $\overset{v}{\to}$). Then $\mathcal{M}(\overline{\R}_0)$ is a metrizable space.
For any $g\in\mathcal{B}_b^+(\overline{\R}_0)$, $\mu\in\mathcal{M}(\overline{\R}_0)$, we write $\mu(g):=\int_{\overline{\R}_0} g(x)\mu(\d x)$.
A sequence of random elements $\nu_n$ in $\mathcal{M}(\overline{\R}_0)$ converges weakly to $\nu$,
denoted as $\nu_n\overset{d}{\to}\nu, $
if and only if for all $g\in C_c^+(\overline{\R}_0)$, $\nu_n(g)$ converges weakly to $\nu(g)$.
Note that, for any $a>0$, $[a,\infty]$ and $[-\infty,-a]$ are compact subsets of $\overline{\R}_0$.

We claim that there exists a non-decreasing function $h_t$ with $h_t\uparrow\infty$ such that
\begin{equation}\label{def-bt}
\lim_{t\to\infty}e^{\lambda t}h_t^{-\alpha}L(h_t)=1.
\end{equation}
In fact,  using \cite[Theorem 1.5.4]{Bingham}, there
exists a non-increasing function $g$ such that
$g(x)\sim x^{-\alpha}L(x)$, as $x\to\infty$. Then $g(x)\to 0$ as $x\to\infty$.
Define
$$h_t:=\inf\{x>0: g(x)\le e^{-\lambda t}\}.$$
It is clear that $h_t$ is non-decreasing and $h_t\uparrow\infty$.
By the definition of $h_t$, one has that, for any $\epsilon>0$, $$g(h_t/(1+\epsilon))\ge  e^{-\lambda t}\ge g(h_t(1+\epsilon)), $$ which implies that
\begin{align*}
(1+\epsilon)^{-\alpha}&=(1+\epsilon)^{-\alpha}\lim_{t\to\infty}\frac{L(h_t)}{L(h_t/(1+\epsilon))}=\lim_{t\to\infty}\frac{g(h_t)}{g(h_t/(1+\epsilon))}\\
&\le \liminf_{t\to\infty}e^{\lambda t} g(h_t)\le\limsup_{t\to\infty}e^{\lambda t} g(h_t)\\
&\le  \lim_{t\to\infty}\frac{g(h_t)}{g(h_t(1+\epsilon))}=(1+\epsilon)^{\alpha} \lim_{t\to\infty}\frac{L(h_t)}{L(h_t(1+\epsilon))}=(1+\epsilon)^{\alpha}.
\end{align*}
Since $\epsilon$ is arbitrary, we get
$$\lim_{t\to\infty}e^{\lambda t}h_t^{-\alpha}L(h_t)=\lim_{t\to\infty}e^{\lambda t} g(h_t)=1.$$
 In particular,
$h_t=e^{\lambda t/\alpha}$ if $L=1.$
In Lemma \ref{lem:vague},  we will prove that
$$
e^{\lambda t}\mP(h_t^{-1}\xi_s\in\cdot)\overset{v}{\to } s
v_\alpha(\cdot),
$$
where
 $$v_\alpha(dx)=q_1 x^{-1-\alpha}{\bf 1}_{(0,\infty)}(x)\d x+q_2 |x|^{-1-\alpha}{\bf 1}_{(-\infty,0)}(x)\d x,$$
with $q_1$ and $q_2$ being nonnegative numbers,  uniquely determined by the following equation:
if $\alpha\neq 1$
$$c_*=\alpha\Gamma(1-\alpha)\left(q_1e^{-i\pi\alpha/2}+q_2e^{i\pi\alpha/2}\right),$$
and if $\alpha=1$
$$q_1=q_2=\Re(c_*)/\pi.$$

Now we are ready to state our main result.  Define
a renormalized version of $\mathbb{X}_t$ by
\begin{equation}\label{def-N}
\mathcal{N}_t:=\sum_{v\in\mathcal{L}_t}\delta_{h_t^{-1}\xi_t^v}.
\end{equation}
In this paper we will investigate
the limit of $\N_t$ as $t\to\infty$.

\begin{theorem} \label{main-theorem} Under $\P$, $\mathcal{N}_t$
converges weakly to a random measure
$\mathcal{N}_\infty\in \mathcal{M}(\overline{\R}_0)$,
defined on some extension $(\Omega, {\cal G}, P)$ of the probability space on which
the branching L\'evy process is defined,
with  Laplace transform given by
\begin{equation}\label{laplace-N}
E(e^{-\N_\infty(g)})
=\E\left(\exp\left\{-W\int_0^\infty e^{-\lambda r}\int_{\mathbb{R}_{0}}\E (1-e^{-Z_{r} g(x)})v_\alpha(dx)\d r\right\}\right),
\quad g\in C_c^+(\overline{\R}_0).
\end{equation}
Moreover,
 $\mathcal{N}_\infty=\sum_{j}T_j\delta_{e_j},$
 where given $W$,
$\sum_{j}\delta_{e_j}$ is a Poisson random measure with intensity $\vartheta Wv_\alpha(dx)$,
$\{T_j,j\ge 1\}$
is a sequence of i.i.d. random variables with common law:
$$
P(T_j=k)
=\vartheta^{-1}\int_0^\infty  e^{-\lambda r}\P(Z_r=k)\d r,\quad k\ge1,$$ where $\vartheta=\int_0^\infty  e^{-\lambda r}\P(Z_r>0)\d r$,
and $\sum_{j}\delta_{e_j}$ and $\{T_j,j\ge 1\}$ are independent.
\end{theorem}

\begin{remark}\label{rek1}
Write $D_f$ for the set of  discontinuity  points of the  function $f$. Then  by Theorem \ref{main-theorem},  we have that $\N_t(f)\overset{d}{\to}\N_\infty(f)$ for any bounded measurable function $f$ on $\overline{\R}_0$ with compact support satisfying $\N_{\infty}(D_f)=0$
$P$-a.s.
Furthermore,  for any $k\ge1$,
$$\left(\N_t(B_1),\N_t(B_2),\cdots,\N_t(B_k)\right)\overset{d}{\to}\left(\N_\infty(B_1),\N_\infty(B_2),\cdots,\N_\infty(B_k)\right),$$
where $\{B_j\}$ are relatively compact subsets of $ \overline{\R}_0$ satisfying   $\N_\infty(\partial B_j)=0$,  $j=1,\cdots,k$,
$P$-a.s.
See \cite[Theorem 4.4]{Kallenberg} for a proof.
\end{remark}

Now we list the positions of all particles alive at time $t$ in a decreasing order:
$$M_{t,1}\ge M_{t,2}\ge \cdots M_{t,Z_t},$$ and for $n>Z_t$, define $M_{t,n}:=-\infty.$
In particular,
$M_{t,1}=\max_{v\in\mathcal{L}_t}\xi^v_t$ is the rightmost position of the particles alive at time $t$. We also order the atoms of $\N_\infty$ as $M_{(1)}\ge M_{(2)}\ge \cdots \ge M_{(k)}\ge \cdots$. Note that on the set $\mathcal{S}$, the number of the atoms of $\N_{\infty}$ is infinite,
and thus
$M_{(k)}, k\geq 1,$ are well defined.
On the set $\mathcal{S}^c$, $\N_{\infty}$ is null, then we define $M_{(k)}=-\infty$ for $k\geq 1$.

Define $\P^*(\cdot):=\P(\cdot|\mathcal{S})$ $(P^*(\cdot):=P(\cdot|\mathcal{S}))$ and let $\E^* (E^*)$ be the corresponding expectation.

\begin{corollary}\label{joint}
For any $n\ge1$,
$$\left(h_t^{-1}M_{t,1},h_t^{-1}M_{t,2},\dots,h_t^{-1}M_{t,n}; \P^*\right)\overset{d}{\to} \left(M_{(1)},M_{(2)},\dots , M_{(n)}; P^*\right).$$
Moreover, $M_{(k)}>0, k\geq 1$, $P^*$-a.s.
\end{corollary}

We write $R_t:=M_{t,1}=\max_{v\in\mathcal{L}_t}\xi^v_t$.

\begin{corollary}\label{cor:rightmost*}
$$
\left(h^{-1}_t\R_t; \P^*\right)
\overset{d}{\to}\left(M_{(1)}; P^*\right),$$
where the law of $(M_{(1)}; P^*)$ is given by
$$
P^*\left(M_{(1)}
\le x\right)=\left\{\begin{array}{ll}
\E^*\left(e^{-\alpha^{-1}q_1\vartheta W x^{-\alpha}}\right), & x>0;\\
0,&x\le 0.
\end{array}\right.$$

\end{corollary}
{\bf Proof:}
Using Corollary \ref{joint}, we get that
$
\left(h^{-1}_t\R_t; \P^*\right)
\overset{d}{\to}\left(M_{(1)}; P^*\right),$ and $M_{(1)}>0$ $P^*$-a.s.
For any $x>0$, we have that
\begin{align*}
P^*\left(M_{(1)}\le x\right)=&P^*\left(\N_\infty(x,\infty)=0\right)=P^*\left(\sum_j {\bf 1}_{(x,\infty)}(e_j)=0\right)\\
=&\E^*\left(e^{-\vartheta W v_\alpha(x,\infty)}\right)
=\E^*\left(e^{-\alpha^{-1}q_1\vartheta W x^{-\alpha}}\right).
\end{align*}
The proof is now complete.
 \hfill$\Box$

\begin{remark}
Similarly, we can order the particles alive at time $t$ in an increasing order:
$L_{t,1}\le L_{t,2}\le \cdots\le L_{t,Z_t}$. Then we can get the corresponding weak convergence of $(L_{t,1}, L_{t,2},\cdots,L_{t,n})$.
\end{remark}

The rest of the paper is organized as follows. In Section 2, we introduce the one large jump principle and give the proof of Theorem \ref{main-theorem} based on Proposition \ref{prop:tildeN1}, which will be proved in
Subsection 2.3.
The proof of Corollary \ref{joint} will be  given  in Section 3.  In Section 4,
we will give more examples satisfying condition (H2)
and conditions which are weaker than (H2),
but sufficient for the main result of this paper.
We will discuss the front position of the Fisher-KPP equation  \eqref{F-KPP} in Section 5.

\section{Proof of  Theorem \ref{main-theorem}}

\subsection{Preliminaries}\label{pre}
Recall that $h_t$ is a function satisfying \eqref{def-bt}.

\begin{lemma}\label{lem:vague}
For any $g\in C_b^0$ and $s>0$,
\begin{equation}\label{vague-conv1}
\lim_{t\to\infty} e^{\lambda t}\mE\left(g(h_t^{-1}\xi_s)\right)=s\int_{\R_0} g(x)v_\alpha(\d x).
\end{equation}
\end{lemma}
{\bf Proof:}
Let $\nu_t$ be the law of $h_t^{-1}\xi_s$. Then by (H2), we have that, as $t\to\infty,$
\begin{equation}\label{e:1}
\exp\left\{e^{\lambda t}\int_\R (e^{i\theta x}-1)\nu_t(\d x)\right\}=\exp\left\{{e^{\lambda t}\left(e^{s\psi(h_t^{-1}\theta)}-1\right)}\right\}\to \exp\left\{s\widetilde{\psi}(\theta)\right\},
\end{equation}
where
$$\widetilde{\psi}(\theta)=\left\{\begin{array}{cc}
-c_*\theta^\alpha, &\theta>0;\\
-\overline{c_*}|\theta|^\alpha,&\theta\le 0.
\end{array}\right.$$
Note that the left side of \eqref{e:1} is the  characteristic
function
of an infinitely divisible random variable $Y_t$ with L\'evy measure $e^{\lambda t}\nu_t$, and by \eqref{stable-exponent}, $e^{s\widetilde{\psi}(\theta)}$ is
the characteristic function of a strictly $\alpha$-stable random variable $Y$ with L\'evy measure $sv_\alpha(dx)$.
Thus $Y_t$ weakly converges to $Y$.
The desired result \eqref{vague-conv1}  follows immediately from \cite[Theorem 8.7 (1)]{Sato}.
\hfill
$\Box$

It is well known (see \cite[Theorem 1.5.6]{Bingham} for instance) that, for any $\epsilon>0$, there
exists $a_\epsilon>0$ such that for any $y>a_{\epsilon}$ and $x>a_{\epsilon}$,
\begin{align}\label{property:L}
 \frac{L(y)}{L(x)}\le (1-\epsilon)^{-1}\max\{(y/x)^{\epsilon}, (y/x)^{-\epsilon}\}.
\end{align}

\begin{lemma}\label{tail-prob}
There exists $c_0>0$ such that for any $s>0$ and $x>2+2a_{0.5}$,
$$G_s(x):={\rm P}(|\xi_s|>x)\le c_0 s x^{-\alpha}L(x).$$
\end{lemma}
{\bf Proof:}
By \cite[(3.3.1)]{Durrett}, we have that, for any $x>2$,
\begin{align*}
\mP(|\xi_s|>x)&\le \frac{x}{2}\int_{-2x^{-1}}^{2x^{-1}}(1-e^{s\psi(\theta)})\,\d \theta\\
&\le s\frac{x}{2}\int_{-2x^{-1}}^{2x^{-1}}\|\psi(\theta)\|\,\d \theta=s\, \int_0^2 \|\psi(\theta/x)\|\,\d \theta,
\end{align*}
where in the last equality, we use the symmetry of $\|\psi(\theta)\|$.
By (H2), it is clear that there exists $c_1>0$ such that
$$\|\psi(\theta)\|\le c_1\theta^{\alpha}L(\theta^{-1}), \quad |\theta|\le 1.$$
Thus, for $x>2+2a_{0.5}$, using \eqref{property:L} with $\epsilon=0.5$, we get that
$$\mP(|\xi_s|>x)\le c_1 s x^{-\alpha}\int_0^2 \theta^{\alpha} L(x/\theta)\,\d\theta\le  2c_1 s x^{-\alpha}L(x)\int_0^2\theta^{\alpha}(\theta^{-1/2}+\theta^{1/2})\,\d\theta.$$
The proof is now complete.
\hfill$\Box$

\begin{remark}\label{rek:tail}
It follows from Lemma \ref{lem:vague} that
\begin{equation}\label{tail-prob1}
\lim_{t\to\infty} e^{\lambda t} {\rm P} (|\xi_s|\ge h_t)= s\frac{q_1+q_2}{\alpha},
\end{equation}
which implies that
\begin{equation}\label{tail-prob2}
{\rm P} (|\xi_s|\ge x)\sim  \frac{q_1+q_2}{\alpha} \, s \, x^{-\alpha}L(x),\quad x\to\infty.
\end{equation}
\end{remark}

Now we recall the many-to-one formula which is  useful in computing  expectations. Here we only list some special  cases which will be used in our paper. See \cite[Theorem 8.5]{Hardy-Harris} for general cases.

\begin{lemma}[Many-to-one formula]\label{many-to-one}
Let $\{n_t\}$ be a Poisson process with parameter $\beta$ on some probability space $(\Omega, {\cal G}, P)$.
Then for any $g\in\mathcal{B}_b^+(\R)$,
$$\E\left(\sum_{v\in\mathcal{L}_t} g(n_t^v)\right)=e^{\lambda t}
E\left(g(n_t)\right),
$$
and for any $0\le s< t$,
$$
\E\left(\sum_{v\in\mathcal{L}_t} {\bf 1}_{b_v\le t-s}\right)=e^{\lambda t}
P\left(n_t-n_{t-s}=0\right)
=e^{\lambda t} e^{-\beta s}.$$
\end{lemma}

\subsection{Proof of the Theorem \ref{main-theorem}}

Recall that,
on some extension $(\Omega, {\cal G}, P)$ of the probability space on which
the branching L\'evy process is defined, given $W$,
$\sum_{j}\delta_{e_j}$ is a Poisson random measure with intensity $\vartheta Wv_\alpha(dx)$,
$\{T_j, j\ge 1\}$ is a sequence of i.i.d. random variables with common law
$$
P(T_j=k)
=\vartheta^{-1}\int_0^\infty  e^{-\lambda r}\P(Z_r=k)\,\d r,\quad k\ge1,$$ where
 $\vartheta=\int_0^\infty  e^{-\lambda r}\P(Z_r>0)\,\d r$,
and $\sum_{j}\delta_{e_j}$
and $\{T_j, j\ge 1\}$ are independent.

\begin{lemma}\label{lem:laplace}
Let $\N_\infty=\sum_{j}T_j\delta_{e_j}$.
Then,
$\N_\infty\in\mathcal{M}(\overline{\R}_0)$ and
 the Laplace transform of $\N_\infty$ is given by
\begin{equation*}
E(e^{-\N_\infty(g)})=\E\left(\exp\left\{-W\int_0^\infty e^{-\lambda r}\int_{\mathbb{R}_{0}}\E (1-e^{-Z_{r} g(x)})v_\alpha(\d x)\,\d r\right\}\right),
\quad g\in C_c^+(\overline{\R}_0).
\end{equation*}
\end{lemma}
{\bf Proof:}
First note that for any $a>0$,  $\vartheta Wv_\alpha([-\infty,-a]\cup [a,\infty])<\infty$. Thus,
given $W$, $\sum_{j}{\bf 1}_{|e_j|\ge a}$ is  Poisson distributed with parameter $\vartheta Wv_\alpha([-\infty,-a]\cup [a,\infty])$, which implies that $\sum_{j}{\bf 1}_{|e_j|\ge a}<\infty$, a.s. Thus by the definition of $\N_\infty$,
 $$
P(\N_\infty([-\infty,-a]\cup [a,\infty])<\infty)=
 P\left(\sum_{j}{\bf 1}_{|e_j|\ge a}<\infty\right)=1.
 $$
 So $\N_\infty\in\mathcal{M}(\overline{\R}_0)$.
Note that
 \begin{align*}
 \phi(\theta)&:=
 E\left(e^{-\theta T_j}\right)
 =\vartheta^{-1}\sum_{k\ge 1}e^{-\theta k}\int_0^\infty  e^{-\lambda r}\P(Z_r=k)\,\d r\\
 &=\vartheta^{-1}\int_0^\infty  e^{-\lambda r}\E\left(e^{-\theta Z_r},Z_r>0\right)\,\d r\\
 &=1-\vartheta^{-1}\int_0^\infty  e^{-\lambda r}\E\left(1-e^{-\theta Z_r}\right)\,\d r.
 \end{align*}
Thus, for any $g\in C_c^+(\overline{\R}_0)$,
\begin{align*}
E\left(e^{-\N_\infty(g)}\right)
&=E\left(e^{-\sum_j T_j g(e_j)}\right)=E\left(\prod_j \phi(g(e_j))\right)\\
&=\E\left(e^{-\vartheta W\int_{\R_0}(1-\phi(g(x)))v_\alpha(\d x)}\right)\\
&=\E\left(\exp\left\{-W\int_0^\infty e^{-\lambda r}\int_{\mathbb{R}_0}\E (1-e^{-Z_{r} g(x)})v_\alpha(\d x)\,\d r\right\}\right).
\end{align*}
The proof is now complete.
\hfill$\Box$

To prove Theorem \ref{main-theorem},
we use the idea  of ``one large jump",  which has been used in
\cite{Durrett83}
and \cite{BHR, BHR2} for branching random walks.
``One large jump" means that  with large probability, for all $v \in \mathcal{L}_t$, at most one of the random variables  $\{|X_{u,t}| : u \in I_v\}$ is bigger than $h_t\theta/t$ ($\theta>0$). Thus to investigate the limit property of $\N_t$,
defined by \eqref{def-N},
 	we will consider another point process:
$$
 \widetilde{\N}_t:
=\sum_{v\in\mathcal{L}_t}\sum_{u\in I_v}\delta_{h_t^{-1}X_{u,t}}.$$

\begin{proposition} \label{prop:tildeN1}
Under $\P$, as $t\to\infty$,
	$$\widetilde{\mathcal{N}}_t\overset{d}{\to}\mathcal{N}_\infty.$$
\end{proposition}

The proof of this proposition is postponed to   the next subsection.
The following lemma formalizes the well-known
one large jump principle (see, e.g., Steps 3 and 4 in Section 2 of \cite{Durrett83}) at
the level of point processes.  Because of Lemma \ref{lemma:tildeN} below,
it is enough to investigate the weak convergence of $\widetilde{\mathcal{N}}_t$, which is much easier compared to that of ${\mathcal{N}}_t$.

  \begin{lemma} \label{lemma:tildeN}
  	Assume $g \in C_c^+(\overline{\R}_0)$. For any $\epsilon>0$,
  	$$\lim_{t\to\infty}\P\left(|\N_t(g)-\widetilde{\N}_t(g)|>\epsilon\right)=0.$$
  \end{lemma}
{\bf Proof:}
Since $g \in C_c^+(\overline{\R}_0)$,
we have Supp$(g) \subset \{x :|x|> \delta\}$ for some $\delta> 0$ .

{\bf Step 1:}
For any $\theta > 0$, let $A_t(\theta)$ denote the event that for all $v \in \mathcal{L}_t$, at most one of the random variables $\{|X_{u,t}| : u \in I_v\}$ is bigger than $h_t\theta/t$. We claim that
\begin{equation}\label{key}
	\P(A_t(\theta)^c)\to0.
\end{equation}
Note that
\begin{align}\label{2.4}
	\P\left(A_t(\theta)^c|\mathcal{F}_t^\mathbb{T}\right)\le \sum_{v\in\mathcal{L}_t}\P\left(\sum_{u\in I_v}{\bf 1}_{\{|X_{u,t}|>h_t\theta/t\}}\ge2|\mathcal{F}_t^\mathbb{T}\right).
\end{align}
By Lemma \ref{tail-prob} and \eqref{property:L} with $\epsilon=0.5$,  we have that for $h_t\theta/t>2+2a_{0.5}$ and $h_t>a_{0.5}$,
\begin{align}\label{3.2}
\P\left(|X_{u,t}|>h_t\theta/t|\mathcal{F}_t^\mathbb{T}\right)&=\mP\left(|\xi_s|>h_t\theta/t\right)|_{s=\tau_{u,t}}\le c_0 \tau_{u,t} h_t^{-\alpha}t^\alpha\theta^{-\alpha}L(h_t\theta/t)\nonumber\\
&\le 2c_0 \theta^{-\alpha} t^{1+\alpha}  h_t^{-\alpha}L(h_t)[(\theta/t)^{1/2}+(\theta/t)^{-1/2}]:=p_t.
\end{align}
Recall that the number of elements in $I_v$ is $n_t^v+1$. Since conditioned on $\mathcal{F}_t^\mathbb{T}$,
$\{X_{u,t}, u\in I_v\}$ are independent, by  \eqref{3.2}, we get  that
\begin{align*}
	\P\left(\sum_{u\in I_v}{\bf 1}_{\{|X_{u,t}|>h_t\theta/t\}}\ge 2|\mathcal{F}_t^\mathbb{T}\right)
	&\le \sum_{m=2}^{n^v_t+1}
	{n^v_t+1\choose m}
	p_t^m=p_t^2\sum_{m=0}^{n^v_t-1}
	{n^v_t+1\choose m+2}
	p_t^m\\
	&\le p_t^2\sum_{m=0}^{n^v_t-1} n_t^v(n^v_t+1)
	{n^v_t-1\choose m}
	p_t^m\\
	&= p_t^2 n_t^v(n^v_t+1)(1+p_t)^{n^v_t-1}.
\end{align*}
Thus  by \eqref{2.4} and the many-to-one formula (Lemma \ref{many-to-one}), we get that
\begin{align}\label{3.3}
\P(A_t(\theta)^c)&=\E\left(\P(A_t(\theta)^c|\mathcal{F}_t^\mathbb{T})\right)\le e^{\lambda t}p_t^2
E\left(n_t(n_t+1)(1+p_t)^{n_t-1}\right)\nonumber\\
&= e^{\lambda t} p_t^2(2\beta+(1+p_t)\beta^2)e^{\beta p_t},
\end{align}
where $n_t$ is a Poisson process with parameter $\beta$ on some probability space $(\Omega, {\cal G}, P)$.
Since $e^{\lambda t}h_t^{-\alpha}L(h_t)\to1$,
\eqref{key} follows from \eqref{3.2} and \eqref{3.3} immediately.

{\bf Step 2:}
Let $\varrho>\beta+1$ to be chosen later.
Let $B_t(\varrho)$ be the event that for all $v \in \mathcal{L}_t$, $n_t^v\le \varrho  t$.
Using the many-to-one formula, we have that
\begin{align*}
	\P(B_t(\varrho)^c)&\le \E(\sum_{v\in\mathcal{L}_t}{\bf 1}_{n_t^v>\varrho  t})=e^{\lambda t}
P(n_t>\varrho t)\le e^{\lambda t}
	\inf_{r>0}e^{-r \varrho t}
E(e^{r n _t})\\
	&=e^{\lambda t} \inf_{r>0}e^{((e^r-1)\beta-r\varrho )t}=e^{\lambda t}e^{-(\varrho(\log \varrho-\log \beta)-\varrho+\beta) t}.
\end{align*}
Choose $\varrho$ large enough so that $\varrho(\log \varrho-\log \beta)-\varrho+\beta>\lambda$, then
$$\lim_{t\to\infty}\P(B_t(\varrho)^c)= 0.$$

{\bf Step 3:} Since $g\in C_c^+(\overline{\R}_0)$, $g$ is uniformly continuous, that is for any $a>0$,
there exists $\eta>0$ such that  $|g(x_1)-g(x_2)|\le a$ whenever $|x_1-x_2|<\eta$.

Now consider $\theta$ small enough such that  $\varrho\theta<\eta\wedge(\delta/2)$.
Let $v'\in I_v$ be such that $|X_{v',t}|=\max_{u\in I_v}\{|X_{u,t}|\}$. We note that, on the event $A_t(\theta)$, $|X_{u,t}|\le \theta h_t/t\le h_t\delta/2$ for any
$u\in I_v\setminus\{v'\}$ and $t>1$,
and thus  $g(X_{u,t}/h_t)=0$, which implies that
$$\widetilde{\N}_t(g)=\sum_{v\in \mathcal{L}_t}\sum_{u\in I_v}g(X_{u,t} /h_t)=\sum_{v\in \mathcal{L}_t}g(X_{v',t}/h_t).$$
Thus it follows that, on the event $A_t(\theta)$,
\begin{align}\label{2.7.1}
	|\N_t(g)-\widetilde{\N}_t(g)|
	&=\left|\sum_{v\in \mathcal{L}_t }\Big[g(\xi^v_t/h_t)-g(X_{v',t}/h_t)\Big]\right|
\end{align}
Since $\xi^v_t=\sum_{u\in I_v}X_{u,t}$, on  the event $A_t(\theta)\cap B_t(\varrho)$,  we have that
$$h_t^{-1}|\xi_t^v-X_{v',t}|
      = h_t^{-1}\left|\sum_{u\in I_v\setminus \{v'\}} X_{u,t}\right|
\le \theta t^{-1} n_t^v\le \varrho \theta<\eta\wedge(\delta/2). $$
Note that if $|X_{v',t}/h_t|\le \delta/2$, then $|\xi_t^v|/h_t< \delta$, which implies that $g(\xi^v_t/h_t)-g(X_{v',t}/h_t)=0$.
Thus we get that
\begin{equation}\label{2.7.2}
	|g(\xi^v_t/h_t)-g(X_{v',t}/h_t)|=|g(\xi^v_t/h_t)-g(X_{v',t}/h_t)|{\bf 1}_{\{|X_{v',t}|>h_t\delta/2\}}\le a{\bf 1}_{\{|X_{v',t}|>h_t\delta/2\}}.
\end{equation}
It follows from \eqref{2.7.1} and \eqref{2.7.2} that, on the event $A_t(\theta)\cap B_t(\varrho)$,
\begin{align*}
	|\N_t(g)-\widetilde{\N}_t(g)|
	&\le a\sum_{v\in \mathcal{L}_t }{\bf 1}_{\{|X_{v',t}| >h_t\delta/2\}}\le   a \sum_{v\in \mathcal{L}_t }\sum_{u\in I_v}{\bf 1}_{\{|X_{u,t}|> h_t\delta/2\}}\\
	&=a \widetilde{\N}_t\left\{[-\infty,-\delta/2)\cup (\delta/2,\infty]\right\}.
\end{align*}
Let $f\in C_c^+(\overline{\R}_0)$
satisfy $f(x)=1$, for $|x|\ge \delta/2$. Then
$$|\N_t(g)-\widetilde{\N}_t(g)|\le a\widetilde{\N}_t(f)).$$

Combining the three steps above, we get that
\begin{align*}
&\limsup_{t\to\infty}\P\left(|\N_t(g)-\widetilde{\N}_t(g)|>\epsilon\right)\\
&\le \limsup_{t\to\infty}\P\left(A_t(\theta)^c\right)+\P\left(B_t(\varrho)^c\right)+\P\left(\widetilde{\N}_t(h)>a^{-1}\epsilon\right)\\
&=\limsup_{t\to\infty}\P\left(\widetilde{\N}_t(f)>a^{-1}\epsilon\right)=
P\left(\N_\infty(f)>a^{-1}\epsilon\right),
\end{align*}
where the final equality
follows from Proposition \ref{prop:tildeN1} (the proof of Proposition \ref{prop:tildeN1} does not use the result in this lemma).
Then letting  $a\to0$, we get the desired result.
\hfill$\Box$

{\bf Proof of Theorem \ref{main-theorem}:} Using Lemma \ref{lem:laplace}, Proposition \ref{prop:tildeN1} and Lemma \ref{lemma:tildeN},
 the results of  Theorem \ref{main-theorem} follow immediately.
\hfill$\Box$

\subsection{Proof of Proposition \ref{prop:tildeN1}}
To prove the weak convergence of $\widetilde{\N}_t$, we first cut the tree at time $t-s$.  We divide the particles born before time $t$ into two parts: the particles born before time $t-s$ and after $t-s$.
Define
$$\widetilde{\N}_{s,t}:=\sum_{v\in\mathcal{L}_t}\sum_{u\in I_v, b_u> t-s}\delta_{h_t^{-1}X_{u,t}}.$$

  \begin{lemma}\label{cut-tree}
  For any $\epsilon>0$ and $g \in C_c^+(\overline{\R}_0)$,
	$$\lim_{s\to\infty}\limsup_{t\to\infty}\P\left(|\widetilde{\N}_t(g)-\widetilde{\N}_{s,t}(g)|>\epsilon\right)=0.$$
\end{lemma}
{\bf Proof:}  Since $g \in C_c^+(\overline{\R}_0)$, We have Supp$(g) \subset \{x :|x|> \delta\}$ for some $\delta> 0$ .

Let $J_{s,t}$ be the event that for all $u$ with $b_u\le t-s$, $|X_{u,t}|\le h_t \delta/2$. On $J_{s,t}$,
$\widetilde{\N}_t(g)-\widetilde{\N}_{s,t}(g)=0$,
 thus  we only need to show that
\begin{align}\label{2.5}
\lim_{s\to\infty}\limsup_{t\to\infty}\P(J_{s.t}^c)=0.
\end{align}
Recall that   $G_s(x):=\mP(|\xi_s|>x)$. By Lemma \ref{tail-prob},
we have that, for $t$ large enough so that  $h_t\delta/2\ge 2+2a_{0.5}$,
\begin{align}\label{est:J}
	\P(J_{s.t}^c)=1-\P(J_{s.t})&=1-\E\left(\prod_{u:b_u\le t-s}(1-G_{\tau_{u,t}}(h_t\delta/2))\right)\nonumber\\
	&\le \E\Big(\sum_{u:b_u\le t-s}G_{\tau_{u,t}}(h_t\delta/2)\Big)\nonumber\\
	&\le c_0h_t^{-\alpha}(\delta/2)^{-\alpha}L(h_t\delta/2)\E\Big(\sum_{u:b_u\le t-s}\tau_{u,t}\Big).
\end{align}
In the first inequality above, we used
the inequality $1-\prod_{i=1}^n((1-x_i))\le \sum_{i=1}^n x_i$, $x_i\in(0,1)$.
By the definition of $\tau_{u,t}$,
\begin{align*}\sum_{u:b_u\le t-s}\tau_{u,t}&=\sum_{u:b_u\le t-s}\int_0^t {\bf 1}_{(b_u,\sigma_u)}(r)\,\d r\\
&=\int_0^{t-s}\sum_{u} {\bf 1}_{(b_u,\sigma_u)}(r)\,\d r+\int_{t-s}^t\sum_{u} {\bf 1}_{b_u<t-s,\sigma_u>r}\,\d r.
\end{align*}
For the first part, noting that $r\in(b_u,\sigma_u)$ is equivalent to $u\in \mathcal{L}_r$,  we get
$$\E\int_0^{t-s}\sum_{u} {\bf 1}_{(b_u,\sigma_u)}(r)\,\d r=\E\int_0^{t-s} Z_r\,\d r=\int_0^{t-s}e^{\lambda r}\,\d r=\lambda^{-1}(e^{\lambda (t-s)}-1).$$
For the second part, using the many-to-one formula, we have that
$$\E\left(\sum_{u} {\bf 1}_{b_u<t-s,\sigma_u>r}\right)=\E\left(\sum_{u\in \mathcal{L}_r} {\bf 1}_{b_u<t-s}\right)=e^{\lambda r} e^{-\beta (r+s-t)}.$$
Thus,
$$\E\int_{t-s}^t\sum_{u} {\bf 1}_{b_u<t-s,\sigma_u>r}\,\d r=\int_{t-s}^t e^{\lambda r}e^{-\beta(r-t+s)}\,\d r=e^{\lambda t}\frac{e^{-\beta s}-e^{-\lambda s}}{\lambda-\beta}.$$
Combining all the above, we get that
\begin{equation}\label{lim:J}
\P(J_{s.t}^c)\le c_0(\delta/2)^{-\alpha}e^{\lambda t}h_t^{-\alpha}L(h_t\delta/2)\Big(\lambda^{-1}e^{-\lambda s}+\frac{e^{-\beta s}-e^{-\lambda s}}{\lambda-\beta}\Big).
\end{equation}
Now first letting $t\to\infty$, and then $s\to\infty$,
we get \eqref{2.5} immediately.
The proof is now complete.
\hfill$\Box$

Now we consider the weak convergence of $\widetilde{\N}_{{s,t}}$.
For $w\in\mathcal{L}_{t-s}$,  let  $\mathcal{L}_t^w:=\{v\in\mathcal{L}_t: w\in I_v\}$ be the set of all the offspring of $w$ at time $t$.
We  rewrite $\widetilde{\N}_{{s,t}}$ as follows:
\begin{align}\label{dec:N}
 \widetilde{\N}_{s,t}=\sum_{w\in\mathcal{L}_{t-s}}M_{s,t}^w,
\end{align}
where
$M_{s,t}^w:=\sum_{v\in\mathcal{L}_t^w}\sum_{u\in I_v,b_u>t-s}\delta_{h_t^{-1}X_{u,t}}$ are i.i.d. with common law
$$M_{s,t}:=
\sum_{v\in\mathcal{L}_s}\sum_{u\in I_v\setminus \{o\}}\delta_{h_t^{-1}X_{u,s}}
=\sum_{u\in D_s }Z^u_{s}\delta_{h_t^{-1}X_{u,s}},$$
where $Z_s^u$ is the number of the offspring of $u$ at time $s$,
and  $D_s=\{u: b_u\le s\}\setminus\{o\}$.

\begin{lemma}\label{lem2}
For any $ j=1,\cdots,n$,
let $\gamma_j(t)$ be a $(0, 1]$-valued function on $(0, \infty)$. Suppose $a_t$ is a positive function with $\lim_{t\to\infty}a_t=\infty$ such that
$\lim_{t\to\infty}a_t(1-\gamma_j(t))\to c_j<\infty.$
Then
$$\lim_{t\to\infty}a_t\left(1-\prod_{j=1}^n \gamma_j(t)\right)\to \sum_{j=1}^n c_j.$$
\end{lemma}
{\bf Proof:}
Note that
$$1-\prod_{j=1}^n \gamma_j(t)=\sum_{j=1}^n\prod_{k=1}^{j-1}\gamma_{k}(t)(1-\gamma_j(t)).$$
Since $\gamma_j(t)\to 1$, thus we get that, as $t\to\infty$,
$$a_t\left(1-\prod_{j=1}^n \gamma_j(t)\right)=\sum_{j=1}^n\prod_{k=1}^{j-1}\gamma_{k}(t) a_t(1-\gamma_j(t))\to \sum_{j=1}^n c_j.$$
\hfill$\Box$

{\bf Proof of Proposition \ref{prop:tildeN1}:}
 By Lemma  \ref{cut-tree}, we only need to consider the convergence of $\widetilde{\N}_{s,t}$. Assume that
 Supp$(g) \subset \{x :|x|> \delta\}$ for some $\delta> 0$.
Using the Markov property and the decomposition of $\widetilde{\N}_{s,t}$ in \eqref{dec:N}, we have that
\begin{align}\label{2.3}
\E\Big(e^{-\widetilde{\N}_{s,t}(g)}\Big)=\E\left([\E(e^{-M_{s,t}(g)})]^{Z_{t-s}}\right).
\end{align}
We claim that
\begin{align}\label{2.2}
\lim_{t\to\infty}	\left(1-\E(e^{-M_{s,t}(g)})\right)e^{\lambda t}
      =\int_{\R_0} \E\left[\sum_{u\in D_s}\tau_{u,s}1-e^{-Z_s^u g(x)}\right]v_\alpha (\d x).
\end{align}
By the definition of $M_{s,t}$, we have that
$$	\left(1-\E(e^{-M_{s,t}(g)}|\mathcal{F}_s^\mathbb{T})\right)e^{\lambda t}=e^{\lambda t}\left(1-\prod_{u\in D_s} \E(e^{-Z_{s}^ug(h_t^{-1}X_{u,s})}|\mathcal{F}_s^\mathbb{T})\right).$$
Note that, given $\mathcal{F}_s^\mathbb{T}$, $X_{u,s}\overset{d}{=}\xi_{\tau_{u,s}}$.
Thus by Lemma  \ref{lem:vague} (with $s$  replaced by $\tau_{u,s}$ and $g$ replaced by $1-e^{-Z^u_s g(x)}$)  we get that
$$e^{\lambda t}\left(1-\E[e^{-Z_s^u g(h_t^{-1}X_{u,s})}|\mathcal{F}_s^\mathbb{T}]\right)\to\tau_{u,s}\int_{\R_0} 1-e^{-Z_s^u g(x)}v_\alpha (\d x),$$
as $t\to\infty.$
Hence it follows from Lemma \ref{lem2} that
\begin{align}\label{3.1}
\lim_{t\to\infty}e^{\lambda t}	\left(1-\E[e^{-M_{s,t}(g)}|\mathcal{F}_s^\mathbb{T}]\right)= \int_{\R_0} \sum_{u\in D_s}\tau_{u,s}[1-e^{-Z_s^u g(x)}]v_\alpha (dx).
\end{align}
Moreover, for $h_t\delta\ge 2+2a_{0.5}$,
\begin{align}\label{3.4}
e^{\lambda t}	\left(1-\E[e^{-M_{s,t}(g)}|\mathcal{F}_s^\mathbb{T}]\right)&\le e^{\lambda t}\E\left(M_{s,t}(g)|\mathcal{F}_s^\mathbb{T}\right)
\le  \|g\|_\infty e^{\lambda t}\sum_{u\in D_s}Z_s^u G_{\tau_{u,s}}(h_t\delta)\nonumber\\
&\le c_0 \|g\|_\infty \delta^{-\alpha}e^{\lambda t}h_t^{-\alpha}L(h_t\delta)\sum_{u\in D_s}\tau_{u,s}Z_s^u\nonumber\\
&\le C\sum_{u\in D_s}\tau_{u,s}Z_s^u,
\end{align}
where $C$ is a constant not depending on $t$. The third inequality follows from Lemma \ref{tail-prob} and the final inequality from the fact $e^{\lambda t}h_t^{-\alpha}L(h_t\delta)\to 1$.
Since $\tau_{u,s}=\int_0^s {\bf 1}_{(b_u,\sigma_u)}(r)\,\d r$,
we have that
\begin{align*}
\E\left(\sum_{u\in D_s}\tau_{u,s}Z_s^u\right)=&
\int_0^s\E\left(\sum_{u\in D_s} {\bf 1}_{(b_u,\sigma_u)}(r) Z^u_s\right)\,\d r
\\ =&\int_0^s\E\left(\sum_{u\in \mathcal{L}_r-\{o\}}  Z^u_s\right)\,\d r\le \int_0^s \E(Z_s)\,\d r
=se^{\lambda s}<\infty.
\end{align*}
Thus by \eqref{3.1}, \eqref{3.4} and the dominated convergence theorem, the claim \eqref{2.2} holds.

By \eqref{2.2} and the fact that $\lim_{t\to\infty}e^{-\lambda t}Z_{t-s}=e^{-\lambda s}W$, we have
$$\lim_{t\to\infty}\left[\E\left(e^{-M_{s,t}(g)}\right)\right]^{Z_{t-s}}=\exp\left\{-e ^{-\lambda s}W\int_{\R_0} \E\Big[\sum_{u\in D_s}\tau_{u,s}(1-e^{-Z_s^u g(x)})\Big]v_\alpha (\d x)\right\}.$$
Thus by \eqref{2.3} and the bounded convergence theorem, we get that
$$\lim_{t\to\infty}
\E\left(e^{-\widetilde{\N}_{s,t}(g)}\right)
=\E\left(\exp\left\{-e ^{-\lambda s}W\int_{\R_0} \E\Big[\sum_{u\in D_s}\tau_{u,s}(1-e^{-Z_s^u g(x)})\Big]v_\alpha (\d x)\right\}\right).$$
By the definition of $\tau_{u,s}$, we have that
\begin{align*}
\sum_{u\in D_s}\tau_{u,s}\left(1-e^{-Z_s^u g(x)}\right)
&=\sum_{u\in D_s}\int_0^s {\bf 1}_{(b_u,\sigma_u)}(r)\,\d r \left(1-e^{-Z_s^u g(x)}\right)
&=\int_0^s \sum_{u\in\mathcal{L}_r\setminus\{o\}}  \left(1-e^{-Z_s^u g(x)}\right)\,\d r.
\end{align*}
Using the Markov property, and the branching property,  $Z^u_s, u\in \mathcal{L}_r$ are i.i.d. with the same distribution as $Z_{s-r}$, and independent with $\mathcal{L}_r$. Thus,
\begin{align}\label{3.6}
\E\sum_{u\in D_s}\tau_{u,s}\left(1-e^{-Z_s^u g(x)}\right)
&=\int_0^s \E\left(Z_r-{\bf 1}_{\{o\in\mathcal{L}_r\}}\right)
\E \left(1-e^{-Z_{s-r} g(x)}\right)\,\d r\nonumber\\
&=\int_0^s \left(e^{\lambda r}-e^{-\beta r}\right)\E \left(1-e^{-Z_{s-r} g(x)}\right)\,\d r,
\end{align}
which implies that
$$e^{-\lambda s}\E\sum_{u\in D_s}\tau_{u,s}\left(1-e^{-Z_s^u g(x)}\right)\to \int_0^\infty e^{-\lambda r}\E \left(1-e^{-Z_{r} g(x)}\right)\,\d r,$$
and
$$e^{-\lambda s}\E\sum_{u\in D_s}\tau_{u,s}\left(1-e^{-Z_s^u g(x)}\right)
\le \int_0^\infty e^{-\lambda r}\E \left(1-e^{-Z_{r} g(x)}\right)\,\d r\le \lambda^{-1} {\bf 1}_{\{|x|>\delta\}}.$$
The final inequality follows from the fact that Supp$(g) \subset \{x :|x|> \delta\}$.
Since $v_\alpha({\bf 1}_{\{|x|>\delta\}})<\infty,$
using the dominated convergence theorem we get that
$$\lim_{s\to\infty}e ^{-\lambda s}\int_{\R_0}\E\Big[\sum_{u\in D_s}\tau_{u,s}( 1-e^{-Z_s^u g(x)})\Big]v_\alpha (\d x)
= \int_0^\infty e^{-\lambda r}\int_{\mathbb{R}_0}\E \left(1-e^{-Z_{r} g(x)}\right)v_\alpha(\d x)\,\d r,$$
which implies that
$$\lim_{s\to\infty}\lim_{t\to\infty}
\E\left(e^{-\widetilde{\N}_{s,t}(g)}\right)
=\E\left(\exp\left\{-W\int_0^\infty e^{-\lambda r}\int_{\mathbb{R}_0}\E \left(1-e^{-Z_{r} g(x)}\right)v_\alpha(\d x)\,\d r\right\}\right).$$
By Lemmas \ref{cut-tree} and  \ref{lem:laplace}, we get that
$$\lim_{t\to\infty}
 \E\left(e^{-\widetilde{\N}_{t}(g)}\right)
=E\left(e^{-\N_\infty(g)}\right).$$
The proof is now complete.
\hfill$\Box$

\section{Joint convergence of the order statistics}

{\bf Proof of Corollary \ref{joint}:}
First, we will show that $M_{(k)}>0$, a.s..
Recall that, given $W$,  $\sum_{j}\delta_{e_j}$ is a Poisson random measure with intensity $\vartheta Wv_\alpha$. Since $v_\alpha(0,\infty)=\infty$,
we have that
$P^*\left(\sum_{j}{\bf 1}_{(0,\infty)}(e_j)=\infty\right)=1$,
which implies that $P^*\left(\N_\infty(0,\infty)=\infty\right)=1$. Thus $M_{(k)}>0$, $P^*$-a.s.

Note that, for any $x\in\overline{\R}_0$, $\N_{\infty}(\{x\})=0,$ a.s..
Since $\{M_{t,k}\le h_t x\}=\{\N_t(x,\infty)\le k-1\}$ for any $x>0$, by Remark \ref{rek1} with $B_k=(x_k,\infty)$, we have that
for any $n\ge 1$  and $x_1,x_2,x_3,\cdots,x_{n}>0$,
\begin{align*}
&\P\left(M_{t,1}\le h_t x_1 ,M_{t,2}\le h_t x_2, M_{t,3}\le h_t x_3,\cdots, M_{t,n}\le  h_t x_{n}\right)\\
=&\P\left(\N_{t}((x_k,\infty)\le k-1, k=1,\cdots,n\right)\\
 \to&P( \N_{\infty}((x_k,\infty))\le k-1,k=1,\cdots,n)\\
=&P(M_{(1)}\le x_1 ,M_{(2)}\le x_2, M_{(3)}\le x_3,\cdots, M_{(n)}\le  x_{n})
\quad \mbox{as } t\to\infty.
\end{align*}
Thus, as $t\to\infty$,
\begin{align}\label{1.3}
&\P^*\left(M_{t,1}\le h_t x_1, M_{t,2}\le h_t x_2,\cdots, M_{t,n}\le  h_t x_{n}\right)\nonumber\\
=&\P(\mathcal{S})^{-1}\left[\P(M_{t,k}\le h_t x_k,k=1,\cdots,n)-\P(M_{t,k}\le h_t x_k,k=1,\cdots,n, \mathcal{S}^c)\right]\nonumber\\
\to &\P(\mathcal{S})^{-1}[P(M_{(k)}\le  x_k,k=1,\cdots,n)-\P( \mathcal{S}^c)]\nonumber\\
=&P^*(M_{(k)}\le  x_k,k=1,\cdots,n),
\end{align}
where in the final equality, we used
the fact that on the event of extinction, $M_{(k)}=-\infty$, $k\ge 1$.

Now we consider the case  $x_1,\cdots,x_n\in \R$ with  $x_i\le 0$ for some $i$, and $x_j>0$, $j\neq i$.
By \eqref{1.3}, we get that,
for any $\epsilon>0$
\begin{align*}
&\limsup_{t\to\infty}\P^*(M_{t,1}\le h_t x_1 ,M_{t,2}\le h_t x_2,\cdots, M_{t,n}\le  h_t x_{n})\\
\le & \lim_{t\to\infty}\P^*(M_{t,j}\le h_t x_j, j\neq i, M_{t,i}\le h_t \epsilon)\\
= &P^*(M_{(j)}\le  x_j,j\neq i, M_{(i)}\le \epsilon).
\end{align*}
The right hand side of the display above tends to 0 as $\epsilon\to0$ since $M_{(i)}>0$ a.s.. Thus
\begin{align}\label{1.4}
\lim_{t\to\infty}\P^*(M_{t,k}\le h_t x_k, k=1,\cdots,n)=0=P^*(M_{(k)}\le  x_k,k=1,\cdots,n).
\end{align}
Similarly, we can get  \eqref{1.4} holds for any $x_1,\cdots,x_n\in \R$.

The proof is now complete.
\hfill$\Box$

\section{Examples and an extension}\label{exam}
In this section,
we first give more examples satisfying  (H2).

\begin{lemma}\label{4.1}
Assume that  $L^*$ is
a positive function on $(0, \infty)$ slowly varying at $\infty$ such that $l_\epsilon(x):=\sup_{y\in(0,x]}y^{\epsilon}L^*(y)<\infty$ for any $\epsilon>0$ and $x>0$.
Then, for any $\epsilon>0$, there exist $c_\epsilon, C_\epsilon>0$ such that for any $y>0$ and $a>c_\epsilon$,
\begin{align}\label{4.2}
\frac{L^*(ay)}{L^*(a)}\le C_\epsilon (y^\epsilon+y^{-\epsilon}).
\end{align}
\end{lemma}
{\bf Proof:}
By \cite[Theorem 1.5.6]{Bingham}, for any $\epsilon>0$,
there exists $c_\epsilon>0$ such that
 for any $a\ge  c_\epsilon$ and $y\ge a^{-1}c_\epsilon$,
\begin{align}\label{4.1.2}
\frac{L^*(ay)}{L^*(a)}\le (1-\epsilon)^{-1}\max\{y^\epsilon, y^{-\epsilon}\}.
\end{align}
Thus for any $a>c_\epsilon$,
\begin{align}\label{4.1.4}
	\frac{L^*(c_\epsilon)}{L^*(a)}\le (1-\epsilon)^{-1}(a/c_\epsilon)^{\epsilon}.
\end{align}
Hence for $a>c_\epsilon$ and $0<y\le a^{-1}c_\epsilon$,  we have that
\begin{align}\label{4.1.3}
\frac{L^*(ay)}{L^*(a)}\le l_\epsilon(c_\epsilon)(ay)^{-\epsilon}/L^*(a)\le \frac{l_\epsilon(c_\epsilon)}{L^*(c_\epsilon)(1-\epsilon)c_\epsilon^\epsilon}y^{-\epsilon}.
\end{align}
Combining \eqref{4.1.2} and \eqref{4.1.3},
there exists $C_\epsilon>0$ such that for any $y>0$ and $a>c_\epsilon$,
$$\frac{L^*(ay)}{L^*(a)}\le C_\epsilon (y^\epsilon+y^{-\epsilon}).$$
\hfill$\Box$

\begin{example} Let
$$n(\d y)=c_1x^{-(1+\alpha)}L^*(x){\bf 1}_{(0,\infty)}(x)\d x+c_2|x|^{-(1+\alpha)}L^*(|x|){\bf 1}_{(-\infty,0)}(x)\d x,$$
where $\alpha\in(0,2), c_1,c_2\ge0$, $c_1+c_2>0$ and $L^*$ is
a positive function on $(0, \infty)$ slowly varying at $\infty$ such that $\sup_{y\in(0,x]}y^{\epsilon}L^*(y)<\infty$ for any $\epsilon>0$ and $x>0$.

(1) For $\alpha\in(0,1)$,
assume that the L\'evy exponent of $\xi$ has the following form:
$$\psi(\theta)=ia\theta-b^2\theta^2+\int (e^{i\theta y}-1)n(\d y),$$
where $a\in\R, b\ge0.$
Using Lemma \ref{4.1} with
$\epsilon\in (0, (1-\alpha)\wedge \alpha)$,
we have that, by the dominated convergence theorem, as $\theta\to 0_+,$
\begin{align*}
\int _0^\infty (e^{i\theta y}-1)n(\d y)&=\theta^\alpha \int _0^\infty (e^{iy}-1)y^{-1-\alpha}L^*(\theta^{-1}y)\,\d y\\
&\sim \theta^\alpha L^*(\theta^{-1})\int _0^\infty (e^{iy}-1)y^{-1-\alpha}\,\d y=-\alpha\Gamma(1-\alpha)e^{-i\pi\alpha/2}\theta^\alpha L^*(\theta^{-1}),
\end{align*}
and
\begin{align*}
\int _{-\infty}^0 (e^{i\theta y}-1)n(\d y)&=\theta^\alpha \int _0^\infty (e^{-iy}-1)y^{-1-\alpha}L^*(\theta^{-1}y)\,\d y\\
&\sim \theta^\alpha L^*(\theta^{-1})\int _0^\infty (e^{-iy}-1)y^{-1-\alpha}\,\d y=-\alpha\Gamma(1-\alpha)e^{i\pi\alpha/2}\theta^\alpha L^*(\theta^{-1}).
\end{align*}
Thus as $\theta\to 0_+,$
$$\psi(\theta)\sim -\alpha\Gamma(1-\alpha)(e^{-i\pi\alpha/2} c_1+e^{i\pi\alpha/2} c_2)\theta^\alpha L^*(\theta^{-1}).$$

(2) For $\alpha\in(1,2)$,
assume that the L\'evy exponent of $\xi$ has the following form:
$$\psi(\theta)=-b^2\theta^2+\int (e^{i\theta y}-1-i\theta y)n(\d y),$$
where $b\ge0.$ Using Lemma \ref{4.1} with
$\epsilon\in (0, (2-\alpha)\wedge(\alpha-1))$,
we have that, by the dominated convergence theorem,  as $\theta\to 0_+,$
\begin{align*}
\int _0^\infty (e^{i\theta y}-1-i\theta y)n(\d y)&=\theta^\alpha \int _0^\infty (e^{i y}-1-iy)y^{-1-\alpha}L^*(\theta^{-1}y)\,\d y\\
&\sim \theta^\alpha L^*(\theta^{-1})\int _0^\infty (e^{iy}-1+iy)y^{-1-\alpha}\,\d y\\
&=-\alpha\Gamma(1-\alpha)e^{-i\pi\alpha/2}\theta^\alpha L^*(\theta^{-1}),
\end{align*}
 and
\begin{align*}
\int _{-\infty}^0 (e^{i\theta y}-1-i\theta y)n(\d y)&=\theta^\alpha \int _0^\infty (e^{-i y}-1+iy)y^{-1-\alpha}L^*(\theta^{-1}y)\,\d y\\
&\sim \theta^\alpha L^*(\theta^{-1})\int _0^\infty (e^{-iy}-1+iy)y^{-1-\alpha}\,\d y\\
&=-\alpha\Gamma(1-\alpha)e^{i\pi\alpha/2}\theta^\alpha L^*(\theta^{-1}).
\end{align*}
Thus as $\theta\to 0_+,$
$$\psi(\theta)\sim -\alpha\Gamma(1-\alpha)(e^{-i\pi\alpha/2} c_1+e^{i\pi\alpha/2} c_2)\theta^\alpha L^*(\theta^{-1}).$$

(3)
For $\alpha=1$, assume that $c_1=c_2$ and the L\'evy exponent of $\xi$ has the following form:
$$\psi(\theta)=ia\theta-b^2\theta^2+\int (e^{i\theta y}-1-i\theta y{\bf 1}_{|y|\le 1})n(\d y),$$
where $a\in\R, b\ge0.$
Since $c_1=c_2$,
we have
\begin{align*}
\int _{-\infty}^\infty (e^{i\theta y}-1-i\theta y{\bf 1}_{|y|\le 1})n(\d y)&=-2c_1\theta \int _0^\infty (1-cos y)y^{-2}L^*(\theta^{-1}y)\,\d y.
\end{align*}
Using Lemma \ref{4.1} with
$\epsilon\in (0, 1)$,
we have that, by the dominated convergence theorem,
$$\lim_{\theta\to 0_+}L^*(\theta^{-1})^{-1}\int _0^\infty (1-cos y)y^{-2}L^*(\theta^{-1}y)\,\d y=\int _0^\infty (1-cos y)y^{-2}\,\d y=\pi/2,$$
which implies that as $\theta\to 0_+,$
$$\psi(\theta)\sim
-(c_1\pi-ia) \theta L^*(\theta^{-1}).
$$
\hfill$\Box$
\end{example}

\begin{remark}[An extension]\label{1-stable}
Checking the proof of Theorem \ref{main-theorem},  we see that
Theorem \ref{main-theorem} holds for  more general branching L\'evy
processes with spatial motions satisfying the following  assumptions:

{\bf (A1)} There
exist a  non-increasing function $h_t$ with $h_t\uparrow \infty$ and
a measure $\pi(\d x)\in \mathcal{M}(\overline{\R}_0)$ such that
$$\lim_{t\to\infty} e^{\lambda t}
\mE(g(h_t^{-1}\xi_s))
=s\int_{\R_0} g(x)\pi(\d x),\quad g\in C_c^+(\overline{\R}_0).$$

{\bf (A2)}
$e^{\lambda t}p_t^2\to0,$ where $p_t:=\sup_{s\le t}\mP(|\xi_s|>h_t\theta/t)$.

{\bf (A3)}
For any  $\theta>0$,
$$\sup_{t>1}\sup_{s\le t}s^{-1}e^{\lambda t}\mP(|\xi_s|>h_t\theta )<\infty.$$

First,  (H2) implies (A1)-(A3).
Next we explain that Theorem \ref{main-theorem} holds under Assumptions (A1)-(A3).
Checking the proof of  Lemma \ref{lemma:tildeN}, we see that Lemma \ref{lemma:tildeN} holds under conditions (A1)-(A3). In fact,
we may replace Lemma \ref{tail-prob} by (A2) to get \eqref{key} (see \eqref{3.2} and \eqref{3.3}).
For the proof of Lemma \ref{cut-tree},
using (A3), we  get that
$$\P(J_{s,t}^c)\le Ce^{-\lambda t}\E \sum_{u:b_u\le t-s}\tau_{u,t},$$
which says that \eqref{est:J} holds.
Thus  \eqref{2.5} holds using the same arguments in Lemma \ref{cut-tree}.
Replacing Lemma \ref{lem:vague} by (A1), we see that Proposition \ref{prop:tildeN1} holds with $v_\alpha$ replaced by $\pi(\d x)$.
So under (A1)-(A3), Theorem \ref{main-theorem} holds with $v_\alpha$ replaced by $\pi(\d x)$.

 An easy example which satisfies (A1)-(A3),
 but does not satisfy (H2) is the non-symmetric 1-stable process.
Assume $\xi$ is a non-symmetric 1-stable process with L\'evy measure
$$n(dx)=c_1x^{-2}{\bf 1}_{(0,\infty)}(x){\rm d} x+c_2|x|^{-2}{\bf 1}_{(-\infty,0)}(x){\rm d} x,$$
where  $c_1,c_2\ge 0$, $c_1+c_2>0,$ and $c_1\neq c_2.$
The L\'evy exponent of $\xi$ is given by, for $\theta>0$
$$\psi(\theta)=-\frac{\pi}{2} (c_1+c_2)\theta-i(c_1-c_2)\theta \log\theta+ia(c_1-c_2)\theta\sim -i(c_1-c_2)\theta \log\theta,\quad \theta\to 0+,$$
where $a$ is constant. Thus $c_*=i(c_1-c_2)$. So $\psi(\theta)$ does not satisfy (H2) since $\Re(c_*)=0$.

By \cite[Section 1.5, Exercise 1]{Bertoin}, we have that
$$\frac{1}{t}\mP(\xi_t\in \cdot)\overset{v}{\to}n(\d x),\quad \mbox{ as }t\to 0.$$
Since $e^{-\lambda t}\xi_s\overset{d}{=}\xi_{s e^{-\lambda t}}+     (c_1-c_2)s\lambda t e^{-\lambda t}$ for $s,t>0$, we have that
$$e^{\lambda t}\mP(e^{-\lambda t}\xi_s\in \cdot)\overset{v}{\to}s \, n(\d x),\quad \mbox{as }t\to\infty.$$
So (A1) holds with $h_t=e^{\lambda t}.$
We claim that, for any $x>0$ and $s>0$,
\begin{equation}\label{claim}
\mP(|\xi_s|>x)\le c(sx^{-1}+s^2x^{-2}+s^2x^{-2}(\log x)^2),
\end{equation}
where $c$ is a constant. Thus it is easy to prove that (A2) and (A3) hold.

In fact, for any $x>0$
\begin{align*}
	\mP(|\xi_s|>x)&\le \frac{x}{2}\int_{-2x^{-1}}^{2x^{-1}}(1-e^{s\psi(\theta)})\,\d \theta=x\int_{0}^{2x^{-1}}(1-\Re(e^{s\psi(\theta)}))\,\d \theta .
\end{align*}
Note that
\begin{align*}
	1-\Re(e^{s\psi(\theta)})&=1-e^{s\Re(\psi(\theta))}\cos[s \Im(\psi(\theta))]\\
	&=1-e^{s\Re(\psi(\theta))}+e^{s\Re(\psi(\theta))}(1-\cos[s \Im(\psi(\theta))])\\
	&\le -s\Re(\psi(\theta))+s^2 [\Im(\psi(\theta))]^2\\
	&=\frac{\pi}{2} (c_1+c_2)s \theta+(c_1-c_2)^2s^2(a- \log\theta)^2\theta^2.
\end{align*}
Thus we have that
\begin{align*}
	\mP(|\xi_s|>x)&\le \pi(c_1+c_2) s\, x^{-1}+(c_1-c_2)^2s^2x^{-2}\int_0^2(a-\log \theta+\log x)^2\theta^2\,\d \theta\\
	&\le  \pi(c_1+c_2) s\, x^{-1}+2(c_1-c_2)^2s^2x^{-2}\int_0^2[(a-\log \theta)^2+(\log x)^2]\theta^2\,\d \theta\\
	&\le c(sx^{-1}+s^2x^{-2}+s^2x^{-2}(\log x)^2),
\end{align*}
which proves the claim \eqref{claim}.
\end{remark}

\section{Front position of Fisher-KPP equation}

Recall that $u_g(t,x)=\E_{\delta_x}\left(e^{-\mathbb{X}_t(g)}\right)
=\E\left(e^{-\sum_{v\in\mathcal{L}_t}g(\xi^v_t+x)}\right)$.
Then $1-u_g(t,x)$ is a mild solution to \eqref{F-KPP}.
For $\theta  \in(0,1)$, the level set $\{x\in \R: 1-u_g(t,x)=\theta\}$ is also called the front of $1-u_g$.
The evolution of the front of
 $1-u_g$ as time goes to $\infty$ is of considerable interest.
Using analytic method,
\cite[Theorem 1.5]{CR}  proved that if the density of $\xi$
 is comparable to that of a symmetric $\alpha$-stable process,
 the front position is exponential in time, which is  in contrast with branching Brownian motion  where it is linear in time. In this paper, we provide a probabilistic proof of
\cite[Theorem 1.5]{CR} using  Corollary \ref{cor:rightmost*},
and also partially generalize it.

\begin{proposition}
\label{5.1}
	\begin{itemize}
		\item [(1)]
Assume that  $a_t$ satisfies  $a_t/h_t\to \infty$ as $t\to\infty$, and that $g$ is a non-negative function  satisfying
		\begin{equation}\label{5.1.1}
e^{\lambda t}\sup_{x\le -a_t/2}g(x)\to 0,
\quad \mbox{ as } t\to\infty.
		\end{equation}
Then
		$$\lim_{t\to\infty}\sup_{x\le -a_t}(1-u_g(t,x))=0. $$
		\item[(2)]
Assume that $c_t$ satisfies  $c_t/h_t\to 0$ as $t\to\infty$, and that $g$ is a non-negative function  satisfying $a_0:=\liminf_{x\to\infty} g(x)>0$.
	Then
				$$\lim_{t\to\infty}\sup_{x\ge -c_t}|u_g(t,x)-\P(\mathcal{S}^c)|=0. $$
	\end{itemize}
\end{proposition}
{\bf Proof:} (1)
Let $g^*(x)=\sup_{y\le-x}g(y)$.
Note that, for $x\leq -a_t$,
\begin{align}\label{5.1.2}
1-u_g(t,x)&=\E\left(1-e^{-\sum_{v\in\mathcal{L}_t}g(\xi^v_t+x)}\right)\nonumber\\
&\le \P(R_t\ge a_t/2) +\E\left(1-e^{-\sum_{v\in\mathcal{L}_t}g(\xi^v_t+x)};R_t<a_t/2\right)\nonumber\\
&\le \P\left(R_t\ge a_t/2) +\E(1-e^{-g^*(a_t/2)Z_t}\right)\nonumber\\
&\le \P(R_t\ge a_t/2)+e^{\lambda t}g^*(a_t/2) ,
\end{align}
where in the second  inequality, we use the fact that, on the event $\{R_t<a_t/2\}$, $\xi^v_t+x<a_t/2-a_t=-a_t/2$ and $g(\xi^v_t+x)\le g^*(a_t/2)$.
By the assumption \eqref{5.1.1},
$e^{\lambda t}g^*(a_t/2)\to 0.$
By Corollary \ref{cor:rightmost*},  one has that
$\P^*(R_t\ge a_t/2)\to0$.  Thus
\begin{align*}
\P(R_t\ge a_t/2)\le \P^*(R_t\ge a_t/2)\P(\mathcal{S})+\P(\|X_t\|>0,\mathcal{S}^c)\to 0,
\end{align*}
as $t\to\infty.$
Thus by \eqref{5.1.2},
$$\lim_{t\to\infty}\sup_{x\le -a_t}(1-u_g(t,x))=0. $$

(2)
Note that
\begin{align}\label{5.1.4}
		|u_g(t,x)-\P(\mathcal{S}^c)|&\le  \E\left(e^{-\sum_{v\in\mathcal{L}_t}g(\xi^v_t+x)};\mathcal{S}\right)+\E\left(1-e^{-\sum_{v\in\mathcal{L}_t}g(\xi^v_t+x)};\mathcal{S}^c\right).
\end{align}
Noticing  that on the event $Z_t=0$, $1-e^{-\sum_{v\in\mathcal{L}_t}g(\xi^v_t+x)}=0$, we get that , for any $x\in\R$
$$\E\left(1-e^{-\sum_{v\in\mathcal{L}_t}g(\xi^v_t+x)};\mathcal{S}^c\right)\le \P(Z_t>0,\mathcal{S}^c)\to 0,$$
as $t\to\infty.$
Let $g_*(x)=\inf_{y\ge x}g(y)$. Since $c_t/h_t\to 0$, for any $\epsilon>0$, there
exists $t_\epsilon>0$ such that $c_t\le \epsilon h_t$ for $t>t_\epsilon$.
For any $t>t_\epsilon$ and $x\geq-c_t$, we have that
\begin{align*}
	\E\left(e^{-\sum_{v\in\mathcal{L}_t}g(\xi^v_t+x)};\mathcal{S}\right)	&\le\E\left(e^{-g_*(c_t)\sum_{v\in\mathcal{L}_t}{\bf 1}_{\xi^v_t>2c_t}};\mathcal{S}\right)\\
	&\le \E\left(e^{-g_*(c_t)\sum_{v\in\mathcal{L}_t}{\bf 1}_{\xi^v_t>2\epsilon h_t}};\mathcal{S}\right)\nonumber\\
	&=\E\left(e^{-g_*(c_t)\N_t(2\epsilon,\infty)};\mathcal{S}\right).
\end{align*}
Thus
\begin{equation}\label{5.1.5}
            \limsup_{t\to\infty}\sup_{x\ge -c_t} |u_g(t,x)-\P(\mathcal{S}^c)|\le
      E\left(e^{-a_0\N_\infty(2\epsilon,\infty)},\mathcal{S}\right).
\end{equation}
Since on the event $\mathcal{S}$,  $\vartheta W v_\alpha(0,\infty)=\infty,$
 thus $\N_\infty(0,\infty)=\infty$.
Now letting $\epsilon\to 0$ in \eqref{5.1.5},  we get the desired result.
\hfill$\Box$

\begin{remark}\label{speed on solution}
Proposition \ref{5.1} is a   slight generalization of  \cite[Theorem 1.5]{CR}.
Assume that $p_0=0$, which ensures that $\P(\mathcal{S}^c)=0$.
If $L=1$,  then  $h_t=e^{\lambda t/\alpha}$,
and we have the following results:
\begin{itemize}
	\item[(1)] Let $g$ be a non-negative measurable function  satisfying
\begin{equation}\label{cond-g}g(x)\le C|x|^{-\alpha},\quad x<0.
\end{equation}
Then for any $\gamma>\lambda/\alpha$
$$
e^{\lambda t} g^*(-e^{\gamma t}/2)\le C2^{\alpha}e^{\lambda t}e^{-\alpha\gamma t}\to 0.
$$
Thus by Proposition \ref{5.1},
we have that
	$$\lim_{t\to\infty}\sup_{x\le -e^{\gamma t}}(1-u_g(t,x))=0. $$
\item[(2)] Assume that $g$ is a non-negative function  satisfying $a_0:=\liminf_{x\to\infty} g(x)>0$.
For any $\gamma<\lambda/\alpha$,
by Proposition \ref{5.1},
we have that
	$$\lim_{t\to\infty}\sup_{x\ge -e^{\gamma t}}u_g(t,x)=0. $$
\end{itemize}
Note that in the notation of \cite{CR}, $\sigma^{**}=\lambda/\alpha$,
and our condition \eqref{cond-g} is equivalent to
\begin{equation}\label{cond-g2}
1-e^{-g(x)}\le C|x|^{-\alpha},\quad x<0,
\end{equation}
for some constant $C$.
If $g$ is nondecreasing,  it is clear that $\liminf_{x\to\infty}g(x)>0$.
Thus when the L\'evy process $\xi$ satisfies (H2) with $L=1$, we can get that the conclusion of \cite[Theorem 1.5]{CR} holds from Proposition \ref{5.1}.
Note that the independent sum of Brownian motion and a symmetric $\alpha$-stable process satisfies (H2) with $L=1$, but its transition density
is not comparable with that of the symmetric $\alpha$-stable process, see \cite{CK2, SZ}.  Note also that the independent sum of a symmetric $\alpha$-stable process
and a symmetric $\beta$-stable process, $0<\alpha<\beta<2$, also satisfies (H2) with $L=1$, but its transition density
is not comparable with that of the symmetric $\alpha$-stable process, see \cite{CK1}.
Note that in this paper we do not need to assume that $g$ is nondecreasing.
Thus Proposition \ref{5.1} partially generalizes \cite[Theorem 1.5]{CR}.
\end{remark}

\vspace{.1in}
\begin{singlespace}

\end{singlespace}

\end{doublespace}

\vskip 0.2truein
\vskip 0.2truein

\noindent{\bf Yan-Xia Ren:} LMAM School of Mathematical Sciences \& Center for
Statistical Science, Peking
University,  Beijing, 100871, P.R. China. Email: {\texttt
yxren@math.pku.edu.cn}

\smallskip
\noindent {\bf Renming Song:} Department of Mathematics,
University of Illinois,
Urbana, IL 61801, U.S.A.
Email: {\texttt rsong@illinois.edu}

\smallskip

\noindent{\bf Rui Zhang:}  School of Mathematical Sciences \& Academy for Multidisciplinary Studies, Capital Normal
University,  Beijing, 100048, P.R. China. Email: {\texttt
zhangrui27@cnu.edu.cn}


\begin{thebibliography}{99}


\bibitem{Aldekon} E. A\"{i}d\'{e}kon. Convergence in law of the minimum of a branching random walk. \emph{Ann. Probab.}, {\bf 41} (2013), 1362--1426.

\bibitem{ABBS} E. A\"{i}d\'{e}kon, J. Berestycki,
\'{E}. Brunet and  Z. Shi.
Branching Brownian motion seen from its tip. \emph{Probab. Theory Related Fields}, {\bf 157} (2013), 405--451.




\bibitem{ABK12} L.-P. Arguin, A. Bovier and N. Kistler. Poissonian statistics in the extremal process of branching Brownian motion. \emph{Ann. Appl. Probab.}, {\bf 22} (2012), 1693--1711.

\bibitem{ABK}L.-P. Arguin, A. Bovier and N. Kistler. The extremal process of branching Brownian motion. \emph{Probab. Theory Related Fields}, {\bf 157} (2013), 535--574.


\bibitem{Athreya-Ney} K. B. Athreya and P. E. Ney. \emph{Branching Processes}. Springer-Verlag Berlin Heidelberg New York, 1972.

\bibitem{Bertoin} J. Bertoin. \emph{L\'evy processes}, Cambridge Univ. Press, Cambridge, 1996.

\bibitem{Bingham} N. H. Bingham, C. M. Goldie and J. L. Teugels. \emph{ Regular Variation}. Cambridge Univ. Press, Cambridge, 1978.

\bibitem{BHR}A. Bhattacharya, R. S. Hazra and P. Roy. Point process convergence for branching random walks with regularly varying steps.
\emph{ Ann. Inst. Henri Poincar\'{e} Probab. Stat.}, {\bf 53} (2017), 802--818.

\bibitem{BHR2} A. Bhattacharya, R. S. Hazra and P. Roy.  Branching random walks, stable point processes and regular variation. \emph{Stochastic Process. Appl.}, {\bf 128} (2018), 182--210.

\bibitem{BMPR}A. Bhattacharya, K. Maulik, Z. Palmowski and P. Roy. Extremes of multi-type branching random walks:  heaviest tail wins.
\emph{ Adv. in Appl. Probab.} {\bf 51} (2019), 514--540.

\bibitem{Bocharov}S. Bocharov. Limiting distribution of particles near the frontier in the catalytic branching Brownian motion.
\emph{Acta Appl. Math.}, {\bf 169} (2020), 433--453.

\bibitem{Bocharov-Harris14}S. Bocharov and S. C. Harris. Branching Brownian motion with catalytic branching at the origin. \emph{Acta Appl. Math.}, {\bf 134} (2014), 201--228.

\bibitem{Bocharov-Harris16}S. Bocharov and S. C. Harris. Limiting distribution of the rightmost particle in catalytic branching Brownian motion. \emph{Electron. Commun. Probab.},
    {\bf 21}, 70 (2016), 12 pp.


\bibitem{Bramson78} M. Bramson.
Maximal displacement of branching Brownian motion. \emph{Comm. Pure Appl. Math.},
{\bf 31} (1978),  531--581.


\bibitem{Bramson} M. Bramson. Convergence of solutions of the Kolmogorov equation to travelling waves.
\emph{Mem. Amer. Math. Soc.},  {\bf 44} (1983), iv+190.



\bibitem{CR}X. Cabr\'{e} and J.-M. Roquejoffre.
The influence of fractional diffusion in Fisher-KPP equations. \emph{Commun. Math. Phys. } {\bf 320} (2013), 679--722.

\bibitem{Carmona-Hu}P. Carmona and Y. Hu. The spread of a catalytic branching random walk.
\emph{Ann. Inst. Henri Poincar\'{e} Probab. Stat.},
{\bf 50} (2014), 327--351.



\bibitem{Chauvin88}B. Chauvin and A. Rouault. KPP equation and supercritical branching Brownian motion in the subcritical speed area. Application to spatial trees.
     \emph{Probab. Theory Related Fields}, {\bf 80} (1988), 299--314.

\bibitem{Chauvin}B. Chauvin and A. Rouault. Supercritical branching Brownian motion and K-P-P equation in the critical speed-area. \emph{Math. Nachr.}, {\bf 149} (1990), 41--59.


\bibitem{CK1} Z.-Q. Chen and T. Kumagai.
Heat kernel estimates for jump processes of mixed type on metric measure spaces.
\emph{Probab. Theory Related Fields}, \textbf{140} (2008), 277--317.

\bibitem{CK2} Z.-Q. Chen and T. Kumagai.
A prior H\"older estimate, parabolic Harnack principle and heat kernel estimates for diffusions with jumps.
\emph{Rev. Mat. Iberoam.}, \textbf{26} (2010), 551--589.


\bibitem{Durrett83} R. Durrett. Maxima of branching random walks. \emph{Z. Wahrsch. Verw. Gebiete}, {\bf 62} (1983), 165--170.

\bibitem{Durrett} R. Durrett. {\it Probability: Theory and Examples.} Fourth edition.  Cambridge Series in Statistical and Probabilistic Mathematics, 31. Cambridge University Press, Cambridge, 2010.



\bibitem{Gantert} N. Gantert. The maximum of a branching random walk with semiexponential increments. \emph{Ann. Probab.}, {\bf 28} (2000), 1219--1229.



\bibitem{Hardy-Harris} R. Hardy and S. C. Harris.
 A spine approach to branching diffusions with applications to Lp-convergence of martingales. \emph{S\'{e}minaire de Probabilit\'{e}s XLII}, 281-330, Lecture Notes in Math., 1979, Springer, Berlin, 2009.

\bibitem{HS}Y. Hu and Z. Shi. Minimal position and critical martingale convergence in branching random walks, and directed polymers on disordered trees.
\emph{Ann.  Probab.}, {\bf 37} (2009), 742--789.




\bibitem{Kallenberg} O. Kallenberg. \emph{Random Measures, Theory and Applications}. Probability Theory and Stochastic Modelling, 77.
 Springer, Cham, 2017.

\bibitem{KPP} A. Kolmogorov, I. Petrovskii and N. Piskounov. \'{E}tude de l'\'{e}quation de la diffusion avec croissance de la quantit\'{e} de la mati$\grave{e}$re at son application $\grave{a}$ un problem biologique. \emph{Moscow Univ. Math. Bull}, {\bf 1} (1937), 1--25.


\bibitem{LS87} S. Lalley and T. Sellke. A conditional limit theorem for the frontier of branching Brownian motion. \emph{Ann. Probab.},
{\bf 15} (1987), 1052--1061.


\bibitem{LS88} S. Lalley and T. Sellke. Travelling waves in inhomogeneous branching Brownian motions I.  \emph{Ann. Probab.},
{\bf 16} (1988), 1051--1062.

\bibitem{LS89} S. Lalley and T. Sellke. Travelling waves in inhomogeneous branching Brownian motions II.  \emph{Ann. Probab.},
{\bf 17} (1989), 116--127.



\bibitem{Madaule} T. Madaule, Convergence in law for the branching random walk seen from its tip. \emph{J. Theoret. Probab.}, {\bf30} (2017), 27--63.



\bibitem{Nishimori-Shiozawa}Y. Nishimori and Y. Shiozawa. Limiting distributions for the maximal displacement of branching Brownian motions. \emph{ J. Math. Soc. Japan},  {\bf 74} (2022), 177--216.


\bibitem{Robert} M. I. Roberts. A simple path to asymptotics for the frontier of a branching Brownian motion. \emph{Ann. Probab.}, {\bf 41} (2013), 3518--3541.

\bibitem{Sato} K.-I. Sato. \emph{L\'evy processes and infinitely divisible distributions}, Cambridge University Press,  2013.

\bibitem{Shiozawa} Y. Shiozawa. Spread rate of branching Brownian motions. \emph{Acta Appl. Math.}, {\bf 155} (2018), 113--150.

\bibitem{Shiozawa2} Y. Shiozawa. Maximal displacement of branching symmetric stable processes.	
 arXiv:2106.15215, 2021.


 \bibitem{SZ} R. Song and Z. Vondracek.
 Parabolic Harnack inequality for the mixture of Brownian motion and stable processes.
 \emph{Tohoku Math. J.}, \textbf{59} (2007), 1--19.



\end{thebibliography}
\end{document}